\documentclass[a4paper,11pt]{article}
\usepackage{paralist}   

\usepackage{amsmath}
\usepackage{amssymb}
\usepackage{amsfonts}
\usepackage{color}
\usepackage{multicol}
\raggedcolumns
\usepackage{verbatim}

\usepackage{graphicx}
\usepackage{subfig}
\usepackage{bm}
\usepackage{soul}

\usepackage{lipsum}
\usepackage{ntheorem}

\renewcommand{\le}{\leqslant}
\renewcommand{\leq}{\leqslant}
\renewcommand{\ge}{\geqslant}
\renewcommand{\geq}{\geqslant}
\newtheorem{theorem}{Theorem}[section]
\renewtheorem*{theorem*}{Theorem}

\newtheorem{corollary}[theorem]{Corollary}
\newtheorem{conjecture}{Conjecture}
\newtheorem{lemma}[theorem]{Lemma}
\renewtheorem*{remark*}{Remark}%[section]

\DeclareMathOperator{\supp}{supp}

\newcommand{\proof}{\textbf{Proof:\ }}
\newcommand{\pbox}{\hfill$\Box$\\}

\newcommand{\R}{\mathbb{R}}
\newcommand{\N}{\mathbb{N}}

\newcommand{\Z}{\mathbb{Z}}

\renewcommand{\S}{\mathbb{S}}

\newcommand{\dist}{\,\mathrm{dist}}

\newcommand{\red}[1]{{\color{red} #1 }}

\newcommand{\rp}{\right)} 
\newcommand{\lp}{\left(}
\newcommand{\rb}{\right]} 
\newcommand{\lb}{\left[}

\newcommand{\bo}{{\cal O}}

\usepackage{geometry}
\usepackage[bookmarks,bookmarksopen=false,% Klappt die Bookmarks in Acrobat aus
                  pdfauthor={Autor},%
                  pdftitle={Titel},%
                  colorlinks=true,%
                  linkcolor=blue,%
                  citecolor=blue,%
                  urlcolor=blue,%
                ]{hyperref}
                

\begin{document}

% The concentration problem for band-limited spherical harmonics
% expansions and the large sieve principle

% A large sieve principle for the concentration problem for
% band-limited spherical harmonics expansions

\title{Concentration estimates for band-limited spherical harmonics
  expansions via the large sieve principle}

\author{ 
M. Speckbacher\thanks{\emph{Univiversit{\'e} de Bordeaux, Institut de
    Math{\'e}matiques de Bordeaux UMR 5251, 351 Cours de la
    Lib{\'e}ration, F-33400, Talence, France, Mail:
    speckbacher@kfs.oeaw.ac.at }}
\ and T. Hrycak\thanks{\emph{Acoustics Research Institute, Austrian
    Academy of Sciences, Wohllebengasse 12-14, 1040 Vienna, Austria}}
}

\date{} 
\maketitle

\begin{abstract}
\noindent We study a concentration problem on the unit sphere $\S^2$
for band-limited spherical harmonics expansions using large sieve
methods. We derive upper bounds for concentration in terms of the
maximum Nyquist density. Our proof uses estimates of the spherical
harmonics coefficients of certain zonal filters. We also demonstrate
an analogue of the classical large sieve inequality for spherical
harmonics expansions.
\end{abstract}

\medskip

\noindent\textbf{MSC2010:} 33C55; 33C45; 46E15; 46E20; 42C10; 11N36
\newline
\textbf{Keywords:} large sieve inequalities, concentration estimates,
spherical harmonics, Legendre polynomials, signal recovery
% Mehler-Heine formula

\section{Introduction}
\subsection{Main contributions}

Let $\S^2$ be the unit sphere in space, $\Omega\subset \S^2$ a
measurable set, and let $\mathcal{S}$ be a Banach subspace of
$L^p(\S^2)$, where $1 < p < \infty$. The concentration problem for the
sphere is concerned with estimating the quantity
\begin{equation} %\label{eq:loc-prob}
\lambda^p_\mathcal{S}(\Omega) := \sup_{f\in
  \mathcal{S}\setminus\{0\}} 
\frac{\int_\Omega|f|^pd\sigma}{\int_{\S^2}|f|^pd\sigma}.
\end{equation}

Following ideas of \cite{dolo92}, we define the \emph{maximum Nyquist
  density} on $\S^2$ as
\begin{equation} \label{def-max-nyqu}
\rho(\Omega,L) = 
\sup_{y\in\S^2} \, \frac{|\Omega\cap 
\mathcal{C}_{t_{L,L}}(y)|}{|\mathcal{C}_{t_{L,L}}(y)|},
\end{equation}
where $t_{L,L}$ denotes the largest zero of the Legendre polynomial
$P_L, L = 1,2,\ldots$, and $\mathcal{C}_{t_{L,L}}(y)$ denotes the
spherical cap with the apex $y \in\S^2$ and the polar angle
$\arccos(t_{L,L})$. A similar concept of density is considered in
\cite{orpri13}.

Let $\mathcal{S}_L$ denote the space of spherical harmonics expansions
with the maximum degree $L$. In this paper, we derive upper bounds for
the concentration constants $\lambda^p_{\mathcal{S}_L}(\Omega)$,
$1<p<\infty$, in terms of the maximum Nyquist density
$\rho(\Omega,L)$. Our approach is to adapt the large sieve principle,
that was first used by Donoho and Logan \cite{dolo92} to study the
concentration problem for band-limited functions on the real line.

% This approach has also been applied recently to study the
% concentration of the short-time Fourier transform on subsets of
% $\R^2$ \cite{abspe17-sampta,abspe18-sieve}.

% concerns the space $\mathcal{S}_L$ of spherical harmonics expansions
% with maximum degree $L$ viewed as a subspace of $L^2(\S^2)$.

Our main result, which is given in Theorem~\ref{thm:l2}, states that
for $L = 1,2, \ldots$
\begin{equation}\label{eq:intro-thm}
\lambda_{\mathcal{S}_L}^2(\Omega)
% = \sup_{f\in\mathcal{S}_L}\frac{\int_\Omega
% |f|^2d\sigma}{\int_{\S^2}|f|^2d\sigma}
\leq B_{L} \cdot\rho(\Omega, L),
\end{equation}
where 
\begin{align}
B_{L}:=(1-t_{L,L})\left(\int_{t_{L,L}}^1 P_L^2(t)dt\right)^{-1}.
% B_{L}:=\frac{1-t_{L,L}}{ \int_{t_{L,L}}^1 P_L^2(t)dt}.
\end{align}
In Lemma~\ref{lemma-binfty}, we show that
\begin{align}
\lim_{L\rightarrow\infty} B_{L} = J_1(j_{0,1})^{-2} \approx
3.71038068570948,
\end{align}
where $J_1$ is the Bessel function of the first kind, and $j_{0,1}$
denotes the smallest positive zero of the Bessel function~$J_0$.
We then derive $L^p$-estimates by interpolation and
duality. Specifically, we demonstrate that for $1<p<\infty$
%
% \begin{equation*}
% \lambda_{\mathcal{S}_L}^p(\Omega)
% \leq \left\{\begin{array}{cl}\hspace{-0.6cm} B_L\cdot\rho(\Omega,L), &
% 2\leq p<\infty\\
% \big(B_L\cdot\rho(\Omega,L)\big)^{p-1},& 1<p<2\end{array}.\right.
% \end{equation*}
%
\begin{equation}
\lambda_{\mathcal{S}_L}^p(\Omega) \leq
\big(B_L\cdot\rho(\Omega,L)\big)^{\min(p-1,\,1)}.
\end{equation}
Donoho and Logan showed that their constants are optimal within their
approach using the Beurling-Selberg function \cite{beu89} and related
extremal functions. Similarly, we show that for $p = 2$, the constant
$B_L$ in \eqref{eq:intro-thm} is also optimal and solves an extremal
problem that can be seen as a spherical analogue of the
Beurling-Selberg problem, and also as a Fourier dual of the problem
considered in \cite[Theorem 4]{liva99}.

>From Theorem~\ref{thm:l2}, we derive an analogue of the classical
large sieve inequality \cite[(2)]{mont78} for spherical harmonics
expansions. Specifically, if
\begin{equation*}
S(x)=\sum_{l=0}^L\sum_{m=-l}^l a_l^m Y_l^m(x),
\end{equation*}
and $x_1,\ldots,x_R\in\S^2$ are $\theta$-separated on the sphere with
$\theta\in(0,\pi]$,
% have a separating angle $\theta\in(0,\pi]$
i.e. $\langle x_k,x_l\rangle \le \cos\theta$, $k\neq l$, then
\begin{equation} %\label{class-ls-intro}
\sum_{k=1}^R|S(x_k)|^2 \leq D(\theta,L) \cdot
\sum_{l=0}^L\sum_{m=-l}^l |a_l^m |^2.
\end{equation}
The constant $D(\theta,L)$ is given explicitly in
Theorem~\ref{thm-class-ls}. Our proof relies on estimating the maximum
number of $\theta$-separated points lying in a spherical cap, which
can be viewed as a packing problem with spherical caps \cite{bamu07}.

\begin{comment}
Donoho and Logan showed that their constants are optimal within their
approach using the Beurling-Selberg function \cite{beu89} and related
extremal functions. Similarly, we show that for $p = 2$, our constant
in \eqref{eq:intro-thm} is also optimal and solves an extremal
problem, that can be seen as a spherical analogue of the
Beurling-Selberg problem and a Fourier dual of the problem considered
in \cite[Theorem 4]{liva99}.

We also show a spherical analogue of \eqref{LS-number} in Theorem
\ref{thm-class-ls}. The proof relies on bounding the maximum number of
$\theta$-separated points in certain spherical caps, which
corresponds to a packing problem of spherical caps \cite{bamu07}.
\end{comment}

\subsection{Previous work}
The concentration problem dealing with the quantity
\begin{equation}\label{eq:conc-slepian}
\mu(\Omega, T):=\sup_{f\in \mathcal{S}_\Omega\setminus\{0\}}
\frac{\int_{-T/2}^{\,\,T/2}\;|f|^2dt}{\int_{\R}\;|f|^2dt},
\end{equation}
where $\mathcal{S}_\Omega = \left\{f\in L^2(\R)\,:\, \widehat{f}(\xi)
= 0, \mbox{ for } |\xi| > \Omega\right\}$, was first studied in a
series of papers by Landau, Slepian and Pollak, now commonly known as
the Bell-Lab papers \cite{lan61, slepo61}.
%
%
% The authors study the eigenvalues and eigenfunctions of the
% concentration operator $P_\mathcal{S}M_{\chi_\Omega}$, where
% $P_\mathcal{S}$ denotes the orthogonal projection on the subspace
% $\mathcal{S}$, and $M_{\chi_\Omega}$ denotes the multiplication by the
% characteristic function of $\Omega$. In this framework, the solution
% to the concentration problem
% %$\lambda^2_{\mathcal{S}}(\Omega)$ 
% is given by the largest eigenvalue of the concentration operator. The
% eigenfunctions of $P_\mathcal{S}M_{\chi_\Omega}$,
%

The largest eigenvalue of the product of the lowpassing operator and
the timelimiting operator corresponds to the solution of
\eqref{eq:conc-slepian}. The eigenfunctions of the product - called
Slepian functions - have appeared in various contexts, for example in
spectral estimation with the multitaper method \cite{tho82,
  abgrro16,abro17}, in time-frequency/time-scale concentration
problems \cite{daub88, daupau88}, and in the study of spatial
concentration of spherical harmonics expansions \cite{sidawi06}. The
Bell-Lab approach has had several generalizations, for example
\cite{abpe15,hogan2012book,hogan2012bandpass,hogan2014}.
%hogan2010sampling, hogan2012book,hogan2012bandpass,hogan2014,hogan2016paley,hogan2016,%hola17

% of functions with compactly supported Hankel transform \cite{abba12},

% Only few explicit solutions of the concentration problem are known,
% and the plethora of contributions following the Bell-Lab approach we
% have mentioned above in no way contradicts this claim. 
There is one common thread throughout the aforementioned papers. They
all exploit specific geometry of concentration domains
% that fit the geometry of the underlying measure space 
in order to solve the concentration problem. For a general
concentration domain, it is hard to explicitly calculate the
eigenvalues following the Bell-Lab theory. Moreover, in many
applications, it is not necessary to know the exact solution to the
concentration problem, and it is enough to have a good estimate. Take
for example the task of reconstructing functions from incomplete
observations.
% Even though the Bell-Lab approach has had many followers, explicit
% solutions of the concentration problem are rare.  Moreover, for many
% applications it is not necessary to know the exact solution to the
% concentration problem, but it suffices to ensure that it is below a
% certain threshold. Reconstructing 
% for example is an ever recurring theme in many signal
% processing problems. 
If a signal is not well-concentrated in a missing region $\Omega$,
then it can be reconstructed by the method of alternating projections,
and the convergence rate is governed by
$\lambda^2_\mathcal{S}(\Omega)<1$, see \cite[Section 4]{dosta89}.

The large sieve principle can be viewed as a class of inequalities
satisfied by trigonometric polynomials $T$ with complex coefficients
\begin{equation*}
T(t)=\sum_{n=1}^N a_n e^{2\pi in t}.
\end{equation*}
Trigonometric polynomials are defined on the interval $[0,1]$
modulo~1, which is endowed with the distance $\dist(t,s) :=
\min_{n\in\Z}|t-s-n|$. If $\delta > 0$ and $t_1, \ldots, t_R \in
    [0,1]$ satisfy
\begin{equation*}
\qquad\dist(t_i, t_j) \ge \delta, \qquad\qquad 1\le i < j \le R,
\end{equation*}
then \cite[Theorem~3]{mont78}
\begin{equation}\label{LS-number}
\sum_{k=1}^R |T(t_k)|^2\leq \lp N-1+\delta^{-1}\rp\sum_{n=1}^N|a_n|^2.
\end{equation}
This is a basic form of the large sieve inequality, and the constant
$N-1+\delta^{-1}$ is sharp. Montgomery \cite{mont78} used
\eqref{LS-number} to study the distribution of prime numbers on large
intervals. A multidimensional version of this estimate can be found in
\cite[Theorem 5]{hova96}.

Donoho and Logan first recognized that \eqref{LS-number} can be used
to \emph{'control the size of trigonometric polynomials on "sparse"
  sets'} \cite{dolo92}, which lead them to derive novel concentration
estimates for band-limited functions. This rationale has recently
inspired a study of the time-frequency concentration problem of the
short-time Fourier transform with Hermite windows
\cite{abspe17-sampta,abspe18-sieve}, and is also a guiding idea for
this contribution.

\section{Preliminaries}
%\label{sec:prel}

Throughout this paper, we use the convention that $x\mbox{ and }y$
denote points on the unit sphere $\S^2$ in space, and $t$ denotes
numbers in the interval $[-1,1]$.

\subsection{Legendre polynomials and the Mehler-Heine formula} 
% We require some results on Legendre polynomials and their
% derivatives. 
%
% We formulate these results and give proofs whenever we could not find
% the exact statement in the literature.
%
\emph{Legendre polynomials} can be defined via the following three term
recurrence \cite[8.914~(1)]{grary07}
\begin{equation}\label{recursion-ordinary}
(n+1)P_{n+1}(t)=(2n+1)tP_n(t)-nP_{n-1}(t), \qquad\qquad n = 1,2,\ldots,
\end{equation}
with $P_0(t)=1$, and $P_1(t)=t$.
%\begin{equation}\label{recursion-deriv}\frac{d}{dx}(P_{n+1}-P_{n-1})
% = (2n+1)P_n\end{equation}
%Moreover, $P_n$ satisfies the following differential equation \cite[8.910 (1)]{grary07}
%\begin{equation}\label{diff-eq}
%(1-t^2)P_n''(t)-2tP'_n(t)+n(n+1)P_n(t)=0.
%\end{equation}
The derivative $P_n'$ satisfies \cite[8.915 (2)]{grary07}
\begin{equation}\label{sum-deriv}
% P_{n}'=(2n-1)P_{n-1}+(2(n-2)-1)P_{n-3}+(2(n-4)-1)P_{n-5}+\ldots .
P_{n}'=(2n-1)P_{n-1}+(2n-5)P_{n-3}+(2n-9)P_{n-5}+\ldots .
\end{equation}
For $t\in[-1,1]$, we have \cite[8.917 (5)]{grary07}
\begin{equation} \label{pn-max}
|P_n(t)| \leq 1,
\end{equation}
which, combined with \eqref{sum-deriv}, gives
\begin{align} \label{pn-prime-max}
|P_n'(t)| \leq (2n-1) + (2n-5) + (2n-9) + \ldots = \frac12\,n(n+1).
\end{align}
%Differentiating \eqref{sum-deriv}, one can show by induction that for every $k = 1, 2, \ldots$, there is a constant $c_k$ such that
%\begin{equation}\label{pn-k-max}
%|P_n^{(k)}(t)|\leq c_k n^{2k}, \qquad\qquad n = 0,1,\ldots.
%\end{equation}

% The values of Legendre polynomials and their derivatives are monotonic
% for certain arguments. 

It is known that all zeros of $P_n$ lie in the interval $(-1,1)$
\cite[18.2(vi)]{NIST10}. For $n\ge 1$, we denote by $t_{n,n}$ the
largest zero of $P_n$. It follows from \cite[18.2(vi)]{NIST10} that
$t_{n,n} < t_{n+1,n+1}$. The following lemma describes monotonicity
properties of Legendre polynomials.

\begin{lemma}\label{ordering-PN}
If $n\ge 1$ and $t\in[t_{n,n},1)$, then for $k=1,\ldots,n$
%\begin{multicols}{2}
\begin{equation}\label{order-i}
%\item
P_k(t)< P_{k-1}(t).
%\item $P_{k}'(t)> P_{k-1}'(t)\geq0$,\label{order-ii}
%\item $P_{k}''(t)> P_{k-1}''(t)\geq0$,\label{order-iii}
%\item $(1-t)P_n'(t)$ is decreasing, \label{order-iv}
%\item $(1-t)P_n''(t)$ is decreasing.\label{order-v}
\end{equation}
Consequently, 
\begin{equation*} %\label{order-ii}
P_k(t) \ge 0.
\end{equation*}
\end{lemma} \proof %Ad $(i)$: 
First we show \eqref{order-i} by induction with respect to~$n$. For
$n=1$, we have $k = 1$, $P_0(t) = 1$, $P_1(t) = t$ and $t_{1,1} = 0$,
so \eqref{order-i} is true. We now assume that \eqref{order-i} holds
for $k = n$. From \eqref{recursion-ordinary}, we have for every $t \in
[t_{n+1,n+1},1)$
\begin{eqnarray*}
(n+1)P_{n+1}(t) &=&(2n+1)t P_n(t)-nP_{n-1}(t)\\
&<& (2n+1) P_n(t)-nP_{n}(t)\\
&=& (n+1)P_n(t).
\end{eqnarray*}
This implies \eqref{order-i} with $k = n+1$ and the inductive proof is
complete.
Since $t_{n,n}$ is the largest zero of $P_n$ and $P_n(1) = 1$, it
follows that $P_n(t) \ge 0$ for $t\in[t_{n,n},1)$. Consequently, for
  $k=1,\ldots,n$, we have
\begin{eqnarray*}
P_k(t) \ge P_n(t) \ge 0.
\end{eqnarray*}

\begin{comment}
\\
Ad $(ii)$: From \cite[8.9.17 and 8.9.18]{NIST10}, we have 
\begin{align*}
&(1-t^2)P_{k}'(t)\ \  \ =\ -k tP_k(t)+kP_{k-1}(t),\\
&(1-t^2)P_{k-1}'(t)=\ k tP_{k-1}(t)-kP_k(t).
\end{align*}
Adding the two equalities and dividing by $1-t^2$ gives
$$
\frac{d}{dt}\big(P_k-P_{k-1}\big)(t) =
\frac{k\big(P_k(t)+P_{k-1}(t)\big)}{t+1}> 0
$$
for every $t\in[t_{n,n},1]$. \\
Ad $(iii)$: Taking the derivative on both sides of \eqref{sum-deriv}
yields using part \eqref{order-ii}
\begin{eqnarray*}
P_{k+1}''(t)&=&(2k+1)P_k'(t)+(2k-3)P_{k-2}'(t)+
(2k-7)P_{k-4}'(t)+\ldots\\
&>&(2k-1)P_{k-1}'(t)+(2k-5)P_{k-3}'(t) +(2k-9)P_{k-5}'(t)+\ldots\\
&=&P_{k}''(t).
\end{eqnarray*}
Ad $(iv)$: Multiplying $((1-t)P_n'(t))'=(1-t)P_n''(t)-P_n'(t)$ by $(1+t)$, we obtain using Legendre's differential equation \eqref{diff-eq}, \eqref{order-i} and \eqref{order-ii} 
\begin{align*}
(1-t^2)P_n''(t)-(1+t)P_n(t)&=2tP_n'(t)-l(l+1)P_n(t)-(1+t)P_n'(t)\\
&=-(1-t)P_n'(t)-l(l+1)P_n(t)\\ &\leq 0.
\end{align*}
%using \eqref{order-i} and \eqref{order-ii}.\\
Ad (v): Use \eqref{sum-deriv} to get
$$
(1-t)P_n''(t)=(1-t)(2n-1)P_{n-1}'(t)+(1-t)(2n-5)P_{n-3}'(t)+\ldots
$$
 Using $t_{1,1}<t_{2,2}<\ldots$, and \eqref{order-iv} thus yield the result.
 \end{comment}
\pbox

% It will later be important for us to evaluate the asymptotics of
% derivatives of the Legendre polynomials evaluated at certain points
% in $[t_{n,n},1]$.

% The largest zero of $P_n$ can be written as
% $t_{n,n}=\cos\theta_{n,1},$ with

For $\theta_{n,1} := \arccos(t_{n,n})$, we have the following
asymptotics \cite[18.16.5]{NIST10}
$$
\theta_{n,1}=\frac{j_{0,1}}{n}+\mathcal{O}(n^{-2}),
$$ 
where $j_{0,1}\approx 2.404825557695772$ denotes the smallest positive
zero of the Bessel function of the first kind $J_0$. Taking the cosine
of both sides yields
\begin{equation}\label{zeros-appr1}
t_{n,n}=1-\frac{j_{0,1}^2}{2n^2}+\mathcal{O}\left(n^{-3}\right).
\end{equation}

The Mehler-Heine formula \cite[18.11.5]{NIST10} describes the
asymptotic behavior of $P_n$ at arguments approaching~$1$
\begin{equation}\label{mehler-heine1}
\lim_{n\rightarrow\infty} P_n \lp1-\frac{z^2}{2n^2}\rp =
% \lim_{n\rightarrow\infty} P_n\Big(\cos\frac{z}{n}\Big) = 
J_0(z).
\end{equation}

\subsection{Spherical harmonics and spherical caps}
Expanding functions in terms of the spherical harmonics is a natural
extension of Fourier series from the unit circle to the three
dimensional sphere.
% \cite{daxu13,sida08}.  
The \emph{spherical harmonics} $Y_l^m $ are given in spherical
coordinates by \cite[14.30.1]{NIST10}
\begin{equation*} %\label{def-spher-harm}
Y_l^m(\theta,\varphi):=\sqrt{\frac{(2l+1)}{4\pi}
\frac{(l-m)!}{(l+m)!}}\,P_l^m(\cos\theta)\,e^{im\varphi},\qquad \theta\in
[0,\pi),\ \varphi\in[0,2\pi),
\end{equation*}
where $0\leq |m|\leq l$, $l = 0, 1, \ldots$, and $P_l^m$ denotes the
associated Legendre polynomial of degree l and order $m$
\cite[14.7.10]{NIST10}
\begin{equation}\label{def-assoc-leg}
P^m_l(t)=\frac{(-1)^{m+l}}{2^l l!}\lp 1-t^2 \rp^{m/2}
\frac{d^{m+l}}{dt^{m+l}}\lp 1-t^2\rp^l.
\end{equation}
In particular, $P_l^0$ coincides with the Legendre polynomial $P_l$
\cite[18.5.5]{NIST10}
\begin{equation*} %\label{def-leg}
P_l(t) = P^0_l(t) =\frac{(-1)^{l}}{2^l l!} \cdot \frac{d^{l}}{dt^{l}}\lp
1-t^2\rp^l.
\end{equation*}
>From \eqref{def-assoc-leg}, we infer that $P^m_l(1) = 0$ if $m\ne 0$.
Consequently, 
\begin{equation*} %\label{assoc-leg1}
P^m_l(1) = \delta_{m}\cdot P^0_l(1) = \delta_{m}\cdot P_l(1) =
\delta_m.
\end{equation*}

% The spherical harmonics are eigenfunctions of the Laplace-Beltrami
% operator.

The family $\{Y_l^m\}_{0\leq|m|\leq l}$ forms an orthonormal basis for
$L^2 (\S^2 )$, where $\S^2$ is equipped with the rotation invariant
surface measure $d\sigma$. The basis coefficients of a function $f\in
L^2 (\S^2 )$ are given by
% $$ \int_{\S^2}Y_l^m
% \overline{Y_k^n}d\sigma=\int_0^\pi\int_0^{2\pi}Y_l^m(\theta,\varphi)
% \overline{Y_k^n(\theta,\varphi)}d\varphi(\theta,\varphi) \sin\theta
% d\varphi d\theta=\delta_{k,l}\delta_{m,n}. $$
%
\begin{equation} \label{def-coeff}
\widehat{f}(l,m)=\int_{\S^2}f(x) \overline{Y_l^m(x)}d\sigma(x) =
\int_0^\pi\int_0^{2\pi}f(\theta,\varphi)
\overline{Y_l^m(\theta,\varphi)}\sin\theta \,d\varphi d\theta.
\end{equation}
In particular, 
\begin{equation}\label{def-spher-harm0}
Y_l^0(\theta,\varphi) = \sqrt{\frac{2l+1}{4\pi}
}\,P_l(\cos\theta),\qquad \theta\in [0,\pi),\ \varphi\in[0,2\pi),
\end{equation}
and
\begin{equation} \label{def-hat0}
\widehat{f}(l,0) = \sqrt{\frac{2l+1}{4\pi}} \int_0^\pi \int_0^{2\pi}
f(\theta,\varphi) \, P_l(\cos\theta) \sin\theta\, d\varphi d\theta.
\end{equation}
Let $\mathcal{S}_L$ be the space of band-limited functions with the
maximum degree $L$, i.e.  $f\in\mathcal{S}_L$, if and only if
$\widehat{f}(l,m)=0$ whenever $l>L$ and $|m|\le l$.

% Note that this space is finite dimensional and therefore invariant
% when being equipped with different $L^p$-norms.

We denote the north pole $(0,0,1)$ of the sphere $\S^2$ by $\eta$.
For $\delta \in[-1,1]$, we define the \emph{spherical cap} with the
apex~$x\in\S^2$ and the polar angle $\arccos\delta$ as follows
\begin{equation*}
\mathcal{C}_\delta(x):=\{y\in\S^2:\ \langle x,y\rangle \ge \delta\}.
\end{equation*}
% The axis of the spherical cap is the line passing through the origin
% and the apex~$x$.
Thus the polar angle is the angle between the ray from the origin to
the apex and the ray from the origin to any point on the boundary of
the cap.
The surface area of the spherical cap $\mathcal{C}_\delta(x)$ does not
depend on the location of the apex~$x$, and is given by the formula
\begin{equation}\label{surf-meas-cap}
\! |\mathcal{C}_\delta(x)| = |\mathcal{C}_\delta(\eta)| =
\int_0^{\pi}\int_0^{2\pi} \chi_{[\delta,1]}(\cos\theta) \sin \theta\,
\,d\varphi d\theta = 2\pi\int_0^{\arccos\delta} \sin \theta \, d\theta =
2\pi(1- \delta).
\end{equation}

% In principle, the $L^2(\S^2)$ concentration problem for functions in
% $\mathcal{S}_L$ is equivalent to the following $(L+1)^2$ dimensional
% quadratic programming problem
% \begin{align*}
% \sup_{f\in \mathcal{S}_L}\frac{\int_\Omega |f|^2d\sigma}{\int_{\S^2}
%   |f|^2d\sigma} = \sup_{a\in\C^{(L+1)^2}}
% \frac{\sum_{l,m}\sum_{k,n}a_l^m \overline{a^n_k}
%   \int_{\Omega}Y_l^m\overline{ Y_k^n} d\sigma}{\|a\|_{\ell^2}^2} =
% \sup_{a\in\C^{(L+1)^2}}\frac{\langle Da,a\rangle}{\|a\|_{\ell^2}^2},
% \end{align*}
% with $D_{lm,kn}:=\int_{\Omega}Y_l^m\overline{ Y_k^n} d\sigma.$ 
% Numerical solutions for general sets $\Omega$ may however require
% substantial effort calculating the entries of the matrix $D$ and is
% often not necessary if appropriate estimates are in place.

\subsection{Convolution on $\S^2$}
% There exist several different notions of convolution on the sphere
% $\S^2$. 

In this paper, we use a concept of convolution with a zonal function
on $\S^2$ that is studied in \cite{kennedy2011}. One advantage of this
approach is that it admits a convolution theorem.

\begin{comment}
$L^2(\S^2)\rightarrow L^2(\S^2)$ with integration over $SO(3)$
$$ f\ast g(x)=\int_{SO(3)}f(s^{-1}x)g(s\eta)ds, $$ 
where $\eta$ denotes the north pole.

$L^2(\S^2)\rightarrow L^2(SO(3))$ (integration over $\S^2$, see filter
banks papers)
$$ f\ast g(s)=\int_{\S^2}f(s y)g(y)d\sigma(y)
$$
\end{comment}
 
Let $g$ be a {\em zonal filter}, i.e. a function on $\S^2 \subset
\R^3$ that only depends on the $z$-coordinate. A zonal filter can be
viewed as a function defined on the interval $[-1,1]$. Thus, with a
slight abuse of notation, we write $g(x)=g(\langle x,\eta\rangle)$,
where $\eta$ denotes the north pole of $\S^2$. 

% In particular, $g \in L^1(\S^2)$ is a zonal filter if and only if
% $\widehat{g}(l,m)=0$ for $0 < |m| \le l$, $l = 1, 2, \ldots$.

We define convolution with the zonal function $g$ as follows
\begin{equation} \label{def-convo}
(f\ast g)(x):=\int_{\S^2}f(y)g(\langle x,y\rangle)d\sigma(y), 
\qquad\qquad x\in\S^2.
\end{equation}
Two numbers $1\le p,q \le \infty$ satisfying $\frac1{p}+\frac1{q} = 1$
are called {\it conjugate exponents}. From H\"older's inequality, we
infer that if $p$ and $q$ are conjugate exponents, then
\begin{align*}
|(f*g)(x)| \le \|f\|_{L^q(\S^2)}\cdot \|g\|_{L^p(\S^2)}, \qquad\qquad
x\in\S^2.
\end{align*}

Since
\begin{equation*}
\|g\|_{L^p(\S^2)} = \|g(\langle \cdot,\eta\rangle)\|_{L^p(\S^2)} =
(2\pi)^\frac1{p}\, \|g\|_{L^p([-1,1])},
\end{equation*}
zonal functions in $L^p(\S^2)$ may be regarded as functions in $L^p\lp
[-1,1]\rp$, $1 \le p\leq \infty$. 

% $\int_{\S^2} |g(\langle x,\eta\rangle)|^p\,d\sigma(x) =
% 2\pi\int_{-1}^1|g(t)|^p\, dt$.
Regarding the Legendre polynomial $P_k$ as a zonal function on $\S^2$,
we have
\begin{equation} \label{pkl0}
\widehat{P_k}(l,0) = \sqrt{\frac{2l+1}{4\pi}} \int_0^{2\pi}
d\varphi  \int_0^{\pi}
P_k(\cos\theta) \, P_l(\cos\theta) \sin\theta\, d\theta
% p= 2\i \delta(k.l) \frac2{2l+1} \sqrt{\frac{2l+1}{4\pi}}
 = \sqrt{\frac{4\pi}{2l+1}}\delta_{k,l}. 
\end{equation}

The following lemma shows that a convolution theorem holds.
\begin{lemma} %\label{lemma-convolution}
If $p$ and $q$ are conjugate exponents, $f\in L^q(\S^2)$ and $g\in
L^p(\S^2)$, then
\begin{equation} \label{conv-thm}
(f\ast g)\ \widehat{}\ (l,m) =
\sqrt{\frac{4\pi}{2l+1}}\, \widehat{f}(l,m)\, \widehat{g}(l,0)
\end{equation}
for $|m| \le l$ and $l = 0, 1, \ldots$.
\end{lemma}

\proof We may assume that $g(x) = P_k(\langle x,\eta\rangle)$, where
$\eta$ is the north pole and $k \ge 0$. The general case follows from
this by a standard approximation argument.
According to an addition theorem for spherical harmonics
\cite[14.30.9]{NIST10}, we have
$$ 
P_k(\langle x, y\rangle) = 
\frac{4\pi}{2k+1}\sum_{n=-k}^k {Y_k^n(x)}\overline{Y_k^n(y)}.
$$ 
Combining this with \eqref{def-coeff} and \eqref{def-convo}, we obtain
\begin{eqnarray*}
(f\ast P_k)\ \widehat{}\ (l,m) &=& \int_{\S^2}\int_{\S^2}f(y)P_k(\langle
  x,y\rangle)d\sigma (y)\;\overline{Y_l^m(x)}d\sigma(x)\\
&=&{\frac{4\pi}{2k+1}}\sum_{n=-k}^k
  \int_{\S^2}f(y)\overline{Y_k^n(y)}d\sigma(y)\int_{\S^2}Y_k^n(x)
  \overline{Y_l^m(x)}d\sigma(x)\\
&=&{\frac{4\pi}{2k+1}}\sum_{n=-k}^k
  \widehat{f}(k,n)\delta_{n,m}\delta_{k,l} = 
{\frac{4\pi}{2k+1}}\widehat{f}(k,m)\delta_{k,l}\\
&=&{\frac{4\pi}{2l+1}}\widehat{f}(l,m)\delta_{k,l}
  =\sqrt{\frac{4\pi}{2l+1}}\widehat{f}(l,m)\widehat{P_k}(l,0).
\end{eqnarray*}
The last equality follows from \eqref{pkl0}.
\pbox

This lemma implies that convolution with a zonal function maps the
space of band-limited functions $S_L$ into itself.

\section{The large sieve inequalities}
% \label{sec:LS-ineq}

\subsection{$\mathbf{L^p}$-bounds for general measures}

Let us denote the space of zonal functions in $ L^p(\S^2)$ that are
supported in the spherical cap $\mathcal{C}_\delta(\eta)$ by
$\mathcal{Z}^p_\delta$. Specifically, for $\delta\in[-1,1]$, we set
$$
\mathcal{Z}^p_\delta:=\big\{g\in L^p(\S^2):\ \supp(g) \subset
        [\delta,1], \ g \mbox{\ is zonal}\big\}.
$$

The following lemma is used in our estimate of
$\lambda^2_{\mathcal{S}_L}(\Omega)$ given in Theorem~\ref{thm:l2}. We
adopt the notation $\|\cdot\|_p = \|\cdot\|_{L^p(\S^2)}$.
\begin{lemma}\label{measure-prop}
Let $\mu$ be a positive $\sigma$-finite measure, and let $1 < p,q <
\infty$ be conjugate exponents. If
$g\in\mathcal{Z}^p_\delta\setminus\{0\}$, then
\begin{equation}\label{sphere-measure}
\int_{\S^2}|f|^pd\mu \leq 
\sup_{h\in \mathcal{S}_L\setminus\{0\}} 
\frac{\|h\|_p^p\|g\|_q^p}{\|h\ast g\|_p^p}
\cdot \|f\|_p^p
\cdot \sup_{y\in \S^2} \mu(\mathcal{C}_{\delta}(y)),
% \int_{\mathcal{C}_{\delta}(y)}d\mu,
\qquad\qquad f\in \mathcal{S}_L.
\end{equation} 
\end{lemma}
\proof We may assume that convolution with $g$ is invertible on
$\mathcal{S}_L$. Otherwise, the first supremum in
\eqref{sphere-measure} is infinite. Since $\supp(g) \subset [\delta,
  1]$, we have
\begin{equation*}
g(\langle x,y\rangle) = g(\langle x,y\rangle) \cdot
\chi_{\mathcal{C}_{\delta}(y)}(x), \qquad\qquad x, y\in\S^2.
\end{equation*}
If $f^\ast\in \mathcal{S}_L$ is a function such that $f=f^\ast \ast
g$, then by H\"older's inequality we have
\begingroup \allowdisplaybreaks
\begin{align} \label{fstar}
\int_{\S^2} |f|^pd\mu&=\int_{\S^2}\left|\int_{\S^2} f^\ast(y)g(\langle
x,y\rangle) \chi_{\mathcal{C}_{\delta}(y)}(x)
d\sigma(y)\right|^pd\mu(x) \nonumber \\
&\leq \int_{\S^2} 
\int_{\S^2}|f^\ast(y)|^p\chi_{\mathcal{C}_{\delta}(y)}(x)d\sigma(y)
\left(\int_{\S^2}|g(\langle x,y\rangle)|^qd\sigma(y)\right)^{p/q}
d\mu(x).
\end{align}
>From rotational invariace of the surface measure~$\sigma$, we infer
that
\begin{equation*} %\label{rotinv}
\left(\int_{\S^2}|g(\langle x,y\rangle)|^qd\sigma(y)\right)^{p/q}
 = \left(\int_{\S^2}|g(\langle \eta,y\rangle)|^qd\sigma(y)\right)^{p/q}
 = \|g\|_q^p, \qquad x\in\S^2.
\end{equation*}
Substituting this into \eqref{fstar} and changing the order of
integration, we obtain
\begin{align*}
\int_{\S^2} |f|^pd\mu
&\le \|g\|_q^p \cdot \int_{\S^2} |f^\ast(y)|^p \;
\mu(\mathcal{C}_{\delta}(y)) \; d\sigma(y)\\
&\leq \|g\|_q^p \cdot \|f^\ast\|^p_p \cdot \sup_{y\in\S^2}
\mu(\mathcal{C}_{\delta}(y))\\
& =  \frac{\|f^\ast\|_p^p\|g\|_q^p}{\|f^\ast\ast
  g\|_p^p} \cdot\|f\|^p_p \cdot \sup_{y\in\S^2}
\mu(\mathcal{C}_{\delta}(y))\\
&\leq \sup_{h\in \mathcal{S}_L\setminus\{0\}}
\frac{\|h\|_p^p\|g\|_q^p}{\|h\ast g\|_p^p} \cdot\|f\|^p_p \cdot
\sup_{y\in\S^2} \mu(\mathcal{C}_{\delta}(y)).
\end{align*}
\endgroup \pbox 

We denote the infimum over $g \in \mathcal{Z}^p_\delta\setminus\{0\}$
of the constants in \eqref{sphere-measure} by
\begin{align} \label{def-cpld}
C_p(L,\delta):=\inf_{g\in \mathcal{Z}^p_\delta\setminus\{0\}} \;
\sup_{h\in
  \mathcal{S}_L\setminus\{0\}}\frac{\|h\|_p^p\|g\|_q^p}{\|h\ast
  g\|_p^p}.
\end{align}
We note that the constant $C_p(L,\delta)$ is the optimal $L^p$-bound
within this approach.

%In the following sections, we consider the two cases $p=2$ and $p=1$.
%, which are of particular interest for \red{signal recovery schemes
% presented in Section ...}.

\subsection{Concentration estimates for
$\bm{\lambda^2_{\mathcal{S}_L}(\Omega)}$}
%%% \subsection{Explicit computation of $\bm{C_2(L,\delta)}$}

In this section, we derive an explicit expression for $C_2(L,\delta)$,
and analyze behavior of this quantity as $L\rightarrow\infty$. In
Theorem~\ref{thm:l2}, we give an upper bound on
$\lambda^2_{\mathcal{S}_L}(\Omega)$ in terms of $C_2(L,\delta)$.

\begin{theorem}\label{thm:l2-const}
If $t_{L,L}\leq \delta< 1$, then the function
$g_\delta:=\chi_{\mathcal{C}_\delta(\eta)}\cdot
P_L\big(\langle\cdot,\eta\rangle\big)$ is a minimizer for the extremal
problem \eqref{def-cpld} defining $C_2(L,\delta)$, and the minimum is
given by
\begin{equation} \label{C2-expl}
C_2(L,\delta) = \left(2\pi\int_{\delta}^1P_L^2(t)dt\right)^{-1}.
\end{equation}
\end{theorem}
\proof First, we simplify the extremal problem \eqref{def-cpld}.  Let
$g\in\mathcal{Z}^2_\delta\setminus\{0\}$. Using the convolution
theorem \eqref{conv-thm} and Parseval's identity, we observe that
\begin{align}
\sup_{h\in \mathcal{S}_L\setminus\{0\}}\frac{\|h\|_2^2\|g\|_2^2}{\|h\ast
  g\|_2^2}&=\sup_{h\in
  \mathcal{S}_L\setminus\{0\}}\|g\|_2^2\|h\|_2^2\lp\sum_{l = 0}^L
\sum_{m = -l}^{l}
\frac{4\pi}{2l+1}|\widehat{h}(l,m)|^2\cdot
|\widehat{g}(l,0)|^2\rp^{-1} \nonumber\\
&=\max_{0\leq l\leq L}
\frac{2l+1}{4\pi}\dfrac{\|g\|_2^2}{|\widehat{g}(l,0)|^2}. \label{extremal2}
\end{align} 
We now show that the constant in \eqref{C2-expl} is attained by the
function $g_\delta$. From \eqref{def-hat0}, we have
$$ 
\sqrt{\frac{4\pi}{2l+1}} \widehat{g_\delta}(l,0) =
2\pi\int_{0}^{\arccos\delta} P_L(\cos\theta)P_l(\cos\theta) \sin\theta
d\theta = 2\pi\int_{\delta}^1P_L(t)P_l(t)dt.  
$$ 
Since $t_{L,L} \le \delta <1$, it follows from Lemma~\ref{ordering-PN}
that
$$ 
\sqrt{\frac{4\pi}{2l+1}} \widehat{g_\delta}(l,0) =
2\pi\int_{\delta}^1P_L(t)P_l(t)dt \ge 2\pi
\int_{\delta}^1P_L^2(t) dx =
\sqrt{\frac{4\pi}{2L+1}}\widehat{g_\delta}(L,0).
$$
Consequently,
\begin{eqnarray*}
\max_{0\leq l\leq L}
\frac{2l+1}{4\pi}\dfrac{\|g_\delta\|_2^2}{|\widehat{g_\delta}(l,0)|^2} 
& = &
\frac{2L+1}{4\pi}\dfrac{\|g_\delta\|_2^2}{|\widehat{g_\delta}(L,0)|^2}
 = 
2\pi\int_{\delta}^1P_L^2(t)dt \cdot
\left(2\pi\int_{\delta}^1P_L^2(t)dt\right)^{-2}\\
& = &
 \left(2\pi\int_{\delta}^1P_L^2(t)dt\right)^{-1}.
\end{eqnarray*}
Finally, we demonstrate that the function $g_\delta$ is a minimizer of
\eqref{extremal2} in $\mathcal{Z}^2_\delta\setminus\{0\}$. From the
Cauchy-Schwarz inequality and \eqref{def-spher-harm0}, we obtain
\begin{align*}
\max_{0\leq l\leq L}
\frac{2l+1}{4\pi}\dfrac{\|g\|_2^2}{|\widehat{g}(l,0)|^2}&\geq
\frac{2L+1}{4\pi}\dfrac{\|g\|_2^2}{|\widehat{g}(L,0)|^2}\geq
\frac{2L+1}{4\pi}\dfrac{\|g\|_2^2}{\|g\|_2^2\cdot
  \|\chi_{\mathcal{C}_\delta(\eta)}\cdot Y_L^0\|_2^2} \\
&=\left(2\pi\int_0^{\arccos\delta}
P_L^2(\cos\theta)\sin\theta d\theta\right)^{-1}
=\left(2\pi\int_{\delta}^1P_L^2(t)dt\right)^{-1}.
\end{align*}
\pbox

We note that a multiple of the function $g_\delta$ is a minimizer of
the following extremal problem: 

find a real valued function $g \in \mathcal{Z}^2_{\delta}$ 
such that $\widehat{g}(l,0)\geq \sqrt{2l+1},\ \ l=0,\ldots,L$, and
whose norm $\|g\|_2$ is minimal.
% \min_{g\in \mathcal{Z}^p_{\delta}}\|g\|^2
% So we need to normalize $g_\delta$ such that $\widehat{g}(L,0)=
% \sqrt{2L+1}$.
>From this perspective, the problem is very similar to
Beurling-Selberg's extremal problem \cite{beu89}, which plays a
central role in the proof of Donoho-Logan's large sieve results for
band-limited functions \cite{dolo92}, and can be seen as a Fourier
side counterpart of an extremal problem considered in \cite[Theorem
  4]{liva99}.

\qquad The following theorem contains our main result.
\begin{theorem}\label{thm:l2} 
Let $\mu$ be a $\sigma$-finite measure, $\Omega\subset\S^2$ be
measurable, and $t_{L,L}\leq\delta<1$. For $L = 1,2,\ldots$ and every
$f\in \mathcal{S}_L$, it holds
\begin{equation}\label{nonuniform-L2-mu}
\int_{\S^2}|f|^2d\mu\leq \left(2\pi\int_\delta^1
P_L^2(t)dt\right)^{-1}\cdot \|f\|_2^2 \cdot
% \int_{\S^2}|f|^2d\sigma \cdot
\sup_{y\in\S^2} \mu(\mathcal{C}_\delta(y)).
\end{equation}
Consequently, 
\begin{equation}\label{uniform-L2}
\lambda^2_{\mathcal{S}_L}(\Omega) \leq B_{L}\cdot \rho(\Omega,L),
\end{equation}
where
\begin{equation}\label{eq-c2l}
B_{L}:=(1-t_{L,L})\left(\int_{t_{L,L}}^1 P_L^2(t)dt\right)^{-1}.
\end{equation}
\end{theorem}
\proof Combining Lemma~\ref{measure-prop} and
Theorem~\ref{thm:l2-const} gives \eqref{nonuniform-L2-mu}. Taking $\mu
= \chi_\Omega d\sigma$ in \eqref{nonuniform-L2-mu} and using
\eqref{surf-meas-cap} and \eqref{def-max-nyqu}, we obtain
\begin{eqnarray*} %\label{uniform-L2proof}
\int_{\Omega}|f|^2d\sigma &\leq&
 \left(2\pi\int_\delta^1 P_L^2(t)dt\right)^{-1} \cdot\|f\|_2^2 \cdot
 \sup_{y\in\S^2} |\Omega\cap \mathcal{C}_\delta(y)|  \\
&\leq&
 \left(2\pi\int_\delta^1 P_L^2(t)dt\right)^{-1} \cdot\|f\|_2^2 \cdot
 \sup_{y\in\S^2} |\Omega\cap \mathcal{C}_{t_{L,L}}(y)| 
\cdot\frac{2\pi(1-t_{L,L})}{| \mathcal{C}_{t_{L,L}}(y)| } \\
& = & (1-t_{L,L}) \left(\int_\delta^1 P_L^2(t)dt\right)^{-1} \cdot \|f\|_2^2
\cdot \rho(\Omega,L),
\end{eqnarray*}
which implies \eqref{uniform-L2}.
\pbox

The behavior of $B_{L}$ for large values of~$L$ is described in the
following lemma.
\begin{lemma} \label{lemma-binfty}
\begin{align} \label{int-pl}
\lim_{L\rightarrow\infty} B_{L} = J_1(j_{0,1})^{-2} \approx
3.71038068570948,
\end{align}
where $J_1$ is the Bessel function of the first kind, and $j_{0,1}$ is
the smallest positive zero of the Bessel function $J_0$.
\end{lemma}
\proof 
We express the integrand in \eqref{eq-c2l} using Taylor's theorem with
the remainder in the Lagrange form
\begin{eqnarray*}
B_{L}^{-1}&=&(1-t_{L,L})^{-1}\int_{t_{L,L}}^1P_L^2(t)dt  \\ &=& \int_0^1
P_L^2\big(1-s(1-t_{L,L})\big)ds  \\
&=& \int_0^1 P_L^2\Big(1-\frac{ j_{0,1}^2}{2L^2}\,s
+h_Ls\Big)ds\\
&=& \int_0^1 \Big[ P_L^2\Big( 1-\frac{ j_{0,1}^2}{2L^2}\,s\Big) +
2h_LsP_L(\xi_s)P_L'(\xi_s) \Big] ds, %\label{int-pl}
\end{eqnarray*}
where $\xi_s\in \Big[1 -\frac{ j_{0,1}^2}{2L^2}s,1 - \frac{
    j_{0,1}^2}{2L^2}s+h_Ls\Big]$, and $h_L=\mathcal{O}(L^{-3})$ in
view of \eqref{zeros-appr1}.
It follows from \eqref{pn-max} and \eqref{pn-prime-max} that
$\|P_L\|_\infty\cdot \|P_L'\|_\infty= \mathcal{O}(L^2)$. From the
Mehler-Heine formula \eqref{mehler-heine1} and the dominated
convergence theorem, we deduce that the integral converges to
$$ 
\int_0^1J_0(j_{0,1}\sqrt{s})^2ds = \frac{2}{j_{0,1}^2} \int_0^{j_{0,1}}s
J_0(s)^2ds =\frac{s^2}{j_{0,1}^2}
\left(J_0(s)^2+J_1(s)^2\right)\Big|_0^{j_{0,1}}= J_1(j_{0,1})^2.
$$ 
The anti-derivative of the function $sJ_0(s)^2$ is given in
\cite[5.54.2]{grary07}. 
\pbox

\subsection{Concentration estimates for
  $\bm{\lambda^p_{\mathcal{S}_L}(\Omega)}$, $\bm{1<p<\infty}$}

Using interpolation and duality arguments, we can extend
\eqref{uniform-L2} to the case ${1<p<\infty}$.
\begin{theorem}
Let $\Omega\subset\S^2$ be measurable and $1<p<\infty$. For
$L=1,2,\ldots$, it holds
\begin{equation*} %\label{eq:other-p}
\lambda_{\mathcal{S}_L}^p(\Omega) =
\sup_{f\in\mathcal{S}_L\setminus\{0\}}
\frac{\int_\Omega|f|^pd\sigma}{\int_{\S^2}|f|^pd\sigma}\leq
\big(B_L\cdot\rho(\Omega,L)\big)^{\min(p-1,1)} .
\end{equation*}
\end{theorem}

\proof The operator
$T_\Omega:\left(\mathcal{S}_L,\|\cdot\|_{L^r(\S^2)}\right)\rightarrow
\left(\mathcal{S}_L,\|\cdot\|_{L^r(\S^2)}\right)$, $T_\Omega
f:=\chi_\Omega\cdot f$, is a contraction for every $1< r < \infty$.
Therefore, the Riesz-Thorin theorem implies that for $2\le p <\infty$
\begin{equation*}
% \lambda_{\mathcal{S}_L}^p(\Omega) = 
\|T_\Omega\|_p \le \|T_\Omega\|_r^{1-\theta} \|T_\Omega\|_2^{\theta}
\le \|T_\Omega\|_2^{\theta},
\end{equation*}
where $r>p$ and $\frac1{p} = \frac{1-\theta}{r} + \frac{\theta}{2}$.
In the limit $r\rightarrow \infty$, we obtain $\|T_\Omega\|_p \le
\|T_\Omega\|_2^{\frac2{p}}$. Consequently,
\begin{equation} \label{eq-thorin1}
\lambda_{\mathcal{S}_L}^p(\Omega) = \|T_\Omega\|^p_p \le
\|T_\Omega\|_2^2 = \lambda_{\mathcal{S}_L}^2(\Omega).
\end{equation}
% 
% As $\lambda_{\mathcal{S}_L}^r(\Omega)\le 1$ for every
% $\Omega\subset \S^2$, $1\le r<\infty$, and $L=1,2,\ldots$, it follows
% from Theorem~\ref{thm:l2} and the Riesz-Thorin theorem that
% $\lambda_{\mathcal{S}_L}^p(\Omega) \le
% \lambda_{\mathcal{S}_L}^2(\Omega) \leq B_L\cdot\rho(\Omega,L)$ for
% every $2\leq p<\infty$.
% 

% The norm of the operator
% $T_\Omega:\left(\mathcal{S}_L,\|\cdot\|_{L^p(\S^2)}\right)\rightarrow
% \left(\mathcal{S}_L,\|\cdot\|_{L^p(\S^2)}\right)$, $T_\Omega
% f:=\chi_\Omega\cdot f$, is thus bounded by
% $\big(B_L\cdot\rho(\Omega,L)\big)^{1/p}$, and so is the norm of the
% 
If $1<p<2$, we consider the 
adjoint operator
$T_\Omega^\ast:\left(\mathcal{S}_L,\|\cdot\|_{L^q(\S^2)}\right)\rightarrow
\left(\mathcal{S}_L,\|\cdot\|_{L^q(\S^2)}\right)$,\ $T_\Omega^\ast
f:=\chi_\Omega\cdot f$, $\frac{1}{p}+\frac{1}{q}=1$. 
Since $2<q<\infty$, we have
\begin{equation}\label{eq-thorin2}
\lambda_{\mathcal{S}_L}^p(\Omega) = \|T_\Omega\|^p_p =
\|T^\ast_\Omega\|^p_q =
\lp\lambda_{\mathcal{S}_L}^q(\Omega)\rp^\frac{p}{q} \le
\lp\lambda_{\mathcal{S}_L}^2(\Omega)\rp^\frac{p}{q}  = 
\lp\lambda_{\mathcal{S}_L}^2(\Omega)\rp^{p-1}.
\end{equation}
The claim now follows from \eqref{eq-thorin1}, \eqref{eq-thorin2} and
\eqref{uniform-L2}. \pbox

\section{The classical large sieve inequality on $\bm{\S^2}$} 
% \label{sec:LS-classic}

In this section, we study the case when the measure $\mu$ in
Theorem~\ref{thm:l2} is a finite sum of Dirac delta distributions,
i.e.  $\mu=\sum_{k=1}^R\delta_{x_k}$. We derive an inequality
analogous to the classical large sieve inequality for trigonometric
polynomials \eqref{LS-number}, see \cite{hova96,mont78}. To this end,
let us assume that the points $x_1, \ldots, x_R$ are
$\theta$-separated on the sphere, i.e. $\langle x_k,x_l\rangle\leq
\cos\theta$, $k\neq l$, for some $\theta \in (0,\pi]$. In other words,
the angle between $x_k$ and $x_l$ is at least~$\theta$. We consider a
spherical harmonics expansion with the maximum degree~$L$
\begin{equation}\label{def-S}
S: = \sum_{l=0}^L\sum_{m=-l}^la_{l}^mY_l^m,
\end{equation}
and intend to find a constant $D = D(\theta, L)$ such that
\begin{equation}\label{pre-constant-ls}
\sum_{k=1}^R|S(x_k)|^2\leq D(\theta,L) \cdot
\sum_{l=0}^L\sum_{m=-l}^l|a_{l}^m|^2.
\end{equation}

%\red{Packing problems like Kepler's conjecture have been studied for
%centuries. On the sphere one of the first packing problems is now
%known as Tammes problem: Fix the number of points on the sphere and
%distribute the points such that the minimum angle between the points
%is maximized \cite{feto99}. Or in other words, distribute the points
%such that one obtains a spherical cap packing with maximal angle. In
%duality one can also fix the angle of the polar cap and ask how many
%polar caps with this angle can be distributed on the sphere
%\textcolor{blue}{find reference}.}

% However, there is no general theory giving good estimates for all
% values of $\alpha$ and $\beta$. 

>From Theorem~\ref{thm:l2}, we obtain the following spherical analogue
of the classical large sieve principle.
\begin{theorem}\label{thm-class-ls}
If $\theta \in (0,\pi]$ and the points $x_1,\ldots,x_R\in\S^2$ are
  $\theta$-separated, then \eqref{pre-constant-ls} holds with the
  constant
\begin{align}\label{classical-ls-const}
D(\theta,L):= \left(2\pi\int_{t_{L,L}}^1P_L^2(t)dt\right)^{-1} \cdot
\frac{1-\cos\frac{\theta}{2} \cdot
  t_{L,L}+\sin\frac{\theta}{2}\cdot\sqrt{1-t_{L,L}^2}}
     {1-\cos\frac{\theta}{2}}.
\end{align}
\end{theorem}
\proof We apply Theorem~\ref{thm:l2} with  $\delta = t_{L,L}$ and $f = S$,
so that
\begin{equation} \label{eq-alm}
\|f\|_2^2 = \|S\|_2^2 = \sum_{l=0}^L\sum_{m=-l}^l|a_l^m|^2.
\end{equation}
It remains to estimate the last factor in \eqref{nonuniform-L2-mu},
that is
\begin{equation} \label{eq-rhodxl}
% \sup_{y\in \S^2}\int_{\mathcal{C}_{t_{L,L}}(y)}d\mu(x) 
\sup_{y\in \S^2}\mu(\mathcal{C}_{t_{L,L}}(y)) = \max_{y\in
  \S^2}\#\{X\cap \mathcal{C}_{t_{L,L}}(y)\},
\end{equation}
where $X:=\{x_k\}_{k=1,\ldots,R}$. Since the points in $X$ are
$\theta$-separated, the angle between every two distinct points in $X$
is at least~$\theta$. Thus the interiors of the spherical caps
$\mathcal{C}_{\cos\frac{\theta}{2}}(x_1), \ldots,
\mathcal{C}_{\cos\frac{\theta}{2}}(x_R)$ with the polar angle
$\frac{\theta}{2}$ are disjoint. Moreover, if $x_k \in
\mathcal{C}_{t_{L,L}}(y)$, then
$\mathcal{C}_{\cos\frac{\theta}{2}}(x_k) \subset
\mathcal{C}_{\cos(\frac{\theta}{2}+\alpha)}(y)$, where $\alpha :=
\arccos(t_{L,L})$.
Therefore, the number of points $x_1, \ldots, x_R$ lying in
$\mathcal{C}_{t_{L,L}}(y)$ does not exceed the maximum number of
spherical caps with the polar angle $\frac{\theta}{2}$ with disjoint
interiors that are contained in a spherical cap with the polar angle
$\frac{\theta}{2}+\alpha$. Comparing the combined areas of the
spherical caps $\mathcal{C}_{\cos\frac{\theta}{2}}(x_1), \ldots,
\mathcal{C}_{\cos\frac{\theta}{2}}(x_R)$ with the area of the
spherical cap $\mathcal{C}_{\cos(\frac{\theta}{2}+\alpha)}(y)$ and
using \eqref{surf-meas-cap}, we obtain
\begin{equation} \label{eq-xctll}
\#\{X\cap \mathcal{C}_{t_{L,L}}(y)\} \le
\frac{|\mathcal{C}_{\cos(\frac{\theta}{2}+\alpha)}(y)|}
     {|\mathcal{C}_{\cos\frac{\theta}{2}}(\cdot)|} =
     \frac{2\pi(1-\cos(\frac{\theta}{2}+\alpha))}
          {2\pi(1-\cos\frac{\theta}{2})}.
\end{equation}
Substituting the following equation 
\begin{equation*}
\cos\Big(\frac{\theta}{2}+\alpha\Big) = \cos\frac{\theta}{2}\cos\alpha
- \sin\frac{\theta}{2}\sin\alpha = \cos\frac{\theta}{2}\cdot t_{L,L} -
\sin\frac{\theta}{2} \cdot \sqrt{1-t_{L,L}^2}
\end{equation*}
into \eqref{eq-xctll}, and taking the maximum over
$y\in\S^2$ yields
\begin{equation} \label{eq-maxys}
\max_{y\in \S^2}\#\{X\cap \mathcal{C}_{t_{L,L}}(y)\} \leq
\frac{1-\cos\frac{\theta}{2} \cdot
  t_{L,L}+\sin\frac{\theta}{2}\cdot\sqrt{1-t_{L,L}^2}}
     {1-\cos\frac{\theta}{2}}.
\end{equation}
Finally, \eqref{pre-constant-ls} follows by combining
\eqref{nonuniform-L2-mu}, \eqref{classical-ls-const}, \eqref{eq-alm},
\eqref{eq-rhodxl} and \eqref{eq-maxys}. \pbox
%
%\indent\qquad In this proof, we deal with a specific packing problem, namely the dual problem to the Tammes problem \cite{tarnai1987}, where the number of spherical caps forming a packing of the sphere is fixed and the minimal polar angle of the spherical caps is maximized \cite{feto99}. Finding bounds on such packing numbers is difficult, see e.g. \cite{bamu07} for estimates of the number of spherical caps contained in another spherical cap in particular cases.

%\red{some comparison to \cite{bamu07}, what are their results, going
%from euclidean to angle distance}

%%%%%%%%%%%%%%%%%%% ASYMPTOTICS
% For $m=0$ we have
% $$
% Y_l^0(\theta,\varphi)=\sqrt{\frac{2l+1}{4\pi}}P_l(\cos\theta)
% $$
% Hence if we choose $a_{l,m}=1$ if $m=0$ and $l=0,1,\ldots,L$ and
% $a_{l,m}=0$ otherwise, then we get
% $$
% S(x)=S(\theta,\varphi)=\sum_{l=0}^L
% \sqrt{\frac{2l+1}{4\pi}}P_l(\cos\theta).
% $$
% If $x=\eta$, the north pole we get by equation (18) that 
% $$
% S(\eta)=\frac{1}{\sqrt{4\pi}}\sum_{l=0}^L \sqrt{2l+1}\asymp
% (L+1)^{\frac{3}{2}}.
% $$
% Consequently,
% $$
% |S(\eta)|^2\asymp (L+1)^{3}=(L+1)^2\sum_{l=0}^L |a_{l,m}|^2.
% $$
% By the last equation in the manuscript, we see that our estimate
% also has the right order with respect to $L$!
%%%%%%%%%%%%%%%%%%%

\qquad We now discuss some basic properties of the expression
appearing in \eqref{classical-ls-const}. From \eqref{int-pl} and
\eqref{zeros-appr1}, we infer that the following quantities are
equivalent up to a constant
\begin{equation} \label{eq-L2}
\Big(2\pi\int_{t_{L,L}}^1P_L^2(t)dt\Big)^{-1}\asymp
(1-t_{L,L})^{-1} \asymp L^2.
\end{equation}
The second factor in \eqref{classical-ls-const} is a decreasing
function of $t_{L,L}$. Since $0 = t_{1,1} \le t_{L,L} <1$, we have
\begin{equation} \label{eq-theta}
1 < \frac{1-\cos\frac{\theta}{2} \cdot
  t_{L,L}+\sin\frac{\theta}{2}\cdot\sqrt{1-t_{L,L}^2}}
{1-\cos\frac{\theta}{2}} \le 
 \frac{1
+\sin\frac{\theta}{2}}
{1-\cos\frac{\theta}{2}}.
\end{equation}

We end this section with a discussion on how close the bound in
Theorem~\ref{thm-class-ls} is to being optimal. 
%
% Following arguments of Montgomery \cite{mont78}, 
%
We derive two elementary lower bounds on the large sieve constants,
and compare them with \eqref{classical-ls-const}. First, let us assume
that we take only one sample~$x_1$ located at the north pole~$\eta$,
and that $a_l^m=\delta_m$, $|m|\le l, l = 0,1,\ldots$. Substituting
\eqref{def-spher-harm0} into \eqref{def-S}, we obtain
\begin{equation*}
S(\eta) = \sum_{l = 0}^L Y_l^0(\eta) = \sum_{l = 0}^L
\sqrt{\frac{2l+1}{4\pi}}.
\end{equation*}
Consequently, the following quantities are equivalent up to a constant
% as $L\rightarrow\infty$, i.e.
\begin{equation} \label{eq-asympL2}
|S(\eta)|^2 \asymp L^3 \asymp L^2 \sum_{l = 0}^L\sum_{m = -l}^l
|a_l^m|^2.
\end{equation}
It follows from \eqref{eq-L2} and \eqref{eq-theta} that $D(\theta,
L)\asymp L^2$ for a fixed $\theta\in(0, \pi]$. Thus \eqref{eq-asympL2}
  implies that for a fixed $\theta$, the bound $D(\theta, L)$ is
  optimal up to a constant factor.

% Thus the constant $D$ in \eqref{pre-constant-ls} satisfies
% \begin{align*}
% D(\theta,L) \geq C_1\cdot L^2
% \end{align*}
% for some constant $C_1>0$ and every $\theta\in[-1,1)$. 
% 
% We thus conclude that for a fixed $\theta$, the constant in
% \eqref{pre-constant-ls} as a function of~$L$ is optimal up to a
% constant factor.\\

% It remains to analyze the behavior of $D$ as a function of $\theta$.
% Let $\{\Phi_k\}_{k=1,\ldots,R}$ be a set of rotations in space such
% that $\langle \Phi_kx,\Phi_m x\rangle \leq \theta$, for every
% $x\in\S^2$ whenever $k\neq m$. From Parseval's identity, we obtain
% \begin{equation}
% \sum_{l=0}^L\sum_{m=-l}^l|a_l^m|^2 = \int_{\S^2}|S(x)|^2d\sigma(x) = 
% \frac1{R}\;\int_{\S^2}\sum_{k=1}^R|S(\Phi_k x)|^2d\sigma(x).
% \end{equation}
% Therefore, there exists $x^\ast\in\S^2$ such that 
% \begin{equation}
% \sum_{k=1}^R|S(\Phi_k x^\ast)|^2 \geq
% \frac{R}{4\pi}\sum_{l=0}^L\sum_{m=-l}^l|a_l^m|^2.
% \end{equation}
% Let $R_{max}(\theta)$ denote the maximum number of $\theta$-separated
% points on $\S^2$. Setting $x_k = \Phi_k x^\ast$, $k = 1, \ldots, R$,
% we see that the constant $D$ in \eqref{pre-constant-ls} satisfies
% \begin{align*}
% D(\theta,L) \geq C_2\cdot R_{max}(\theta)
% \end{align*}
% for some constant $C_2>0$. It holds
% 
% 

% \tr{Given $\theta$-separated points $x_1, \ldots, x_R$, how does one
%   construct rotations $\Phi_1, \ldots, \Phi_R$ such that $\langle
%   \Phi_kx,\Phi_m x\rangle \leq \theta$, for every $x\in\S^2$? }\\ 

It remains to analyze the behavior of $D(\theta, L)$ as a function of
$\theta$ for a fixed~$L$. Let $R_{max}(\theta)$ denote the maximum
number of $\theta$-separated points on $\S^2$. It is known
\cite[p. 121]{gamal1987}, \cite[(24)]{wyner1965} that
\begin{equation} \label{rmax}
R_{max}(\theta) \ge \frac2{1-\cos\theta}.
\end{equation}
For a fixed $\theta$, let $x_1, \ldots, x_{R_{max}(\theta)}\in\S^2$ be
$\theta$-separated, and $a_l^m = 0$, $|m|\le l, l = 0,1,\ldots$,
except for $a_0^0 = 1$. According to \eqref{def-spher-harm0}, we have
%%% $S(x_1) = \ldots = S(x_R) = \frac1{\sqrt{4\pi}}$.
\begin{equation} \label{rmax4pi}
\sum_{k=1}^R |S(x_k)|^2 = \frac{R_{max}(\theta)}{{4\pi}}.
\end{equation}
It follows from \eqref{classical-ls-const} that $D(\theta, L)\asymp
\frac1{1-\cos\theta}$ for a fixed~$L$. Thus from \eqref{rmax} and
\eqref{rmax4pi}, we conclude that also for a fixed~$L$, the bound
$D(\theta, L)$ is within a constant factor from being optimal.

The inequality \eqref{rmax} has a simple proof. If the points
$x_1,\ldots,x_{R_{max}(\theta)}$ on $\S^2$ are $\theta$-separated,
then the union of the spherical caps $\mathcal{C}_{\cos\theta}(x_1)$,
$\ldots$, $\mathcal{C}_{\cos\theta}(x_{R_{max}(\theta)})$ covers the
unit sphere. Otherwise, one could find an additional point on $\S^2$
that is $\theta$-separated from the points
$x_1,\ldots,x_{R_{max}(\theta)}$. Comparing the areas of the caps with
that of the unit sphere, we obtain
\begin{equation} \label{rmaxproof}
R_{max}(\theta)\cdot 2\pi (1-\cos\theta)\geq 4\pi,
\end{equation}
which is equivalent to \eqref{rmax}.

\section*{Acknowledgements}

The authors would like to thank Lu{\'i}s Daniel Abreu for posing the
problem and for his valuable comments and suggestions. We are also
grateful to Uju{\'e} Etayo for a reference about the spherical cap
packing problem.

T.H. was supported by the Innovationsfonds 'Forschung, Wissenschaft
und Gesellschaft' of the Austrian Academy of Sciences on the project
``Railway vibrations from tunnels''.

M.S. was supported by the Austrian Science Fund (FWF) through the
START-project FLAME ('Frames and Linear Operators for Acoustical
Modeling and Parameter Estimation'; Y 551-N13) and  an Erwin-Schr{\"{o}}dinger Fellowship (J-4254).

\bibliographystyle{plain}
\bibliography{paperbib1}
 
\end{document}